\documentclass[a4paper]{article}
\usepackage{latexsym,amssymb,amsmath,amsthm}
\usepackage[all]{xy}
\newcommand{\Spec}{\text{\it Spec}}
\newcommand{\Q}{{\mathbb Q}}
\newcommand{\Z}{{\mathbb Z}}
\newcommand{\ab}{\text{\it ab}}
\newcommand{\et}{\text{\it et}}
\newcommand{\fl}{\text{\it fl}\,}
\newcommand{\cd}{\text{\it cd}\:}
\newcommand{\scd}{\text{\it scd}\:}

\newtheorem{theorem}{Theorem}[section]
\newtheorem{corollary}[theorem]{Corollary}
\newtheorem{lemma}[theorem]{Lemma}
\newtheorem{proposition}[theorem]{Proposition}
\newtheorem{definition}[theorem]{Definition}
\def\freeproduct{\mathop{\lower.4mm\hbox{\emar \symbol{3}}}\limits}
 
 \def\freeproductmed{\mathop{\lower.2mm\hbox{\emas \symbol{3}}}\limits}
 \def\ressum{\mathop{\hbox{${\displaystyle\bigoplus}'$}}\limits}

\font\emar = cmsy10 scaled\magstep4 \font\emas = cmsy10 scaled\magstep2

\def\ressumsmall{\mathop{\hbox{${\bigoplus}'$}}}

\title{\bf Circular sets of prime numbers and $p$-extensions of the rationals}
\author{by Alexander Schmidt}
\date{}

\begin{document} 
\maketitle

\begin{minipage}{11cm}
\footnotesize {\it Abstract:} Let $p$ be an odd prime number and let $S$ be a finite set of prime numbers congruent to $1$ modulo $p$. We prove that the group $G_S(\Q)(p)$ has cohomological dimension~$2$ if the linking diagram attached to $S$ and $p$ satisfies a certain technical condition, and we show that $G_S(\Q)(p)$ is a duality group in these cases.  Furthermore, we investigate the decomposition behaviour of primes in the extension $\Q_S(p)/\Q$ and we relate the cohomology of $G_S(\Q)(p)$ to the {\'e}tale cohomology of the scheme $\Spec(\Z)-S$. Finally, we calculate the dualizing module.
\end{minipage}

\section{Introduction}
Let $k$ be a number field, $p$ \/a prime number and $S$ a finite set of places of $k$. The pro-$p$-group $G_S(k)(p)=G(k_S(p)/k)$, i.e.\ the Galois group of the maximal $p$-extension of $k$ which is unramified outside $S$, contains valuable information on the arithmetic of the number field $k$. If all places dividing $p$ are in $S$, then we have some structural knowledge on $G_S(k)(p)$, in particular, it is of cohomological dimension less or equal to $2$ (if $p=2$ one has to require that $S$ contains no real place, \cite{Sc3}),  and it is often a so-called duality group, see \cite{NSW}, X, \S 7. Furthermore, the cohomology of $G_S(k)(p)$ coincides with the {\'e}tale cohomology of the arithmetic curve $\Spec({\cal O}_k) -S$ in this case.

\medskip
In the opposite case, when $S$ contains no prime dividing $p$, only little is known. By a famous theorem of Golod and \v{S}afarevi\v{c}, $G_S(k)(p)$ may be infinite. A conjecture due to Fontaine and Mazur \cite{FM} asserts that $G_S(k)(p)$ has no infinite quotient which is an analytic pro-$p$-group. So far, nothing was known on the cohomological dimension of $G_S(k)(p)$ and on the relation between its cohomology and the {\'e}tale cohomology of the scheme $\Spec({\cal O}_k) -S$. 

\medskip
Recently, J.~Labute \cite{La}  showed that pro-$p$-groups with a certain kind of relation structure have cohomological dimension $2$.  By a result of H. Koch \cite{Ko}, $G_S(\Q)(p)$ has such a relation structure if the set of prime numbers $S$ satisfies a certain technical condition. In this way, Labute obtained first examples of pairs $(p,S)$ with $p\notin S$ and $\cd G_S(\Q)(p)=2$, e.g.\ $p=3$, $S=\{7,19,61, 163\} $.

\medskip
The objective of this paper is to use arithmetic methods in order to extend Labute's  result. First of all, we weaken the condition on $S$ which implies cohomological dimension~$2$ (and strict cohomological dimension~$3$!) and we show that $G_S(\Q)(p)$ is a duality group in these cases.  Furthermore, we investigate the decomposition behaviour of primes in the extension $\Q_S(p)/\Q$ and we relate the cohomology of $G_S(\Q)(p)$ to the {\'e}tale cohomology of the scheme $\Spec(\Z)-S$. Finally, we calculate the dualizing module.

\section{Statement of results}

Let $p$ be an odd prime number, $S$ a finite set of prime numbers not containing $p$ and $G_S(p)=G_S(\Q)(p)$ the Galois group of the maximal $p$-extension $\Q_S(p)$ of $\Q$ which is unramified outside $S$. Besides $p$, only prime numbers congruent to $1$ modulo $p$ can ramify in a $p$-extension of\/ $\Q$, and we assume that all primes in $S$ have this property. Then $G_S(p)$ is a pro-$p$-group with $n$ generators and $n$ relations, where $n=\#S$ (see lemma~\ref{relrank}).  

\medskip
Inspired by some analogies between knots and prime numbers (cf.\ \cite{Mo}),
J.~Labute \cite{La} introduced the notion of the linking diagram $\Gamma (S)(p)$  attached to $p$ and $S$ and showed that $\cd G_S(p)=2$ if $\Gamma (S)(p)$ is a `non-singular circuit'. Roughly speaking, this means that there is an ordering $S=\{q_1,q_2,\ldots,q_n\}$ such that $q_1q_2\cdots q_nq_1$ is a circuit in $\Gamma(S)(p)$ (plus two technical conditions, see section~\ref{lab} for the definition).

\medskip
We generalize Labute's result by showing 

\begin{theorem} \label{2} Let $p$ be an odd prime number and let $S$ be a finite set of prime numbers congruent to $1$ modulo $p$. Assume there exists a subset $T\subset S$ such that the following conditions are satisfied.
\begin{enumerate}
\item[\rm (i)] $\Gamma(T)(p)$ is a non-singular circuit.
\item[\rm (ii)] For each $q \in S \backslash T$ there exists a directed path in $\Gamma(S)(p)$ starting in $q$ and ending with a prime in $T$.
\end{enumerate}
Then $\cd G_S(p)=2$.
\end{theorem}

\noindent
{\it Remarks.} 1. Condition (ii) of Theorem~\ref{2} can be weakened, see section~\ref{lab}.

\smallskip\noindent
2. Given $p$, one can construct examples of sets $S$ of arbitrary cardinality $\#S \geq 4$ with $\cd G_S(p)=2$ .

\bigskip\noindent
{\it Example.} 
For  $p=3$ and $S=\{7,13,19,61,163\}$, the linking diagram has the following shape

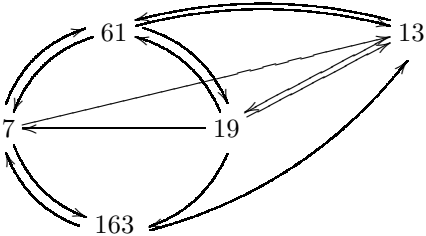
\[
\xymatrix{&61\ar@/^2ex/[rrr]\ar@/_2ex/[dl] \ar@<1ex>@/^2ex/[dr]&&&13\ar@/_2ex/@<-.5ex>[lll] \ar@<-.5ex>[dll]\\
7\ar@<1ex>@/^2ex/[ur]\ar@/_2ex/[dr]\ar[urrrr]&&19\ar[ll]\ar@/_2ex/[ul]\ar@<1ex>@/^2ex/[dl]\ar[urr]\\
&163\ar@<1ex>@/^2ex/[ul]\ar@<-1.3ex>@/_3ex/[uurrr]
}
\]
The linking diagram associated to the subset $T=\{7,19,61,163\}$ is a non-singular circuit, and we obtain $\cd G_S(3)=2$ in this case.

\bigskip
The proof of Theorem~\ref{2} uses arithmetic properties of $G_S(p)$ in order to enlarge the set of prime numbers $S$ without changing the cohomological dimension of $G_S(p)$. In particular, we show

\pagebreak

\begin{theorem}\label{1} Let $p$ be an odd prime number and let $S$ be a finite set of prime numbers congruent to $1$ modulo $p$. Assume that $G_S(p)\neq 1$ and $\cd G_S(p)\leq 2$. Then the following holds.

\begin{enumerate}
\item[\rm (i)] $\cd G_S(p)= 2$ and $\scd G_S(p) = 3$.
\item[\rm (ii)]$G_S(p)$ is a pro-$p$ duality group (of dimension $2$).
\item[\rm (iii)] For all $\ell \in S$, $\Q_S(p)$ realizes the maximal $p$-extension of\/ $\Q_\ell$, i.e.\ (after choosing a prime  above $\ell$ in $\bar \Q$), the image of the natural inclusion $\Q_S(p)\hookrightarrow \Q_\ell (p)$ is dense.
\item[\rm (iv)] The scheme $X=\Spec(\Z)-S$ is a $K(\pi,1)$ for $p$ and the {\'e}tale topology, i.e.\ for any $p$-primary $G_S(p)$-module $M$, considered as a locally constant {\'e}tale sheaf on $X$, the natural homomorphism
\[
H^i (G_S(p), M) \rightarrow H^i_{\et} (X,M)
\]
is an isomorphism for all $i$.
\end{enumerate}
\end{theorem}

\noindent
{\it Remarks.} 1. If $S$ consists of a single prime number, then $G_S(p)$ is finite, hence $\#S\geq 2$ is necessary for the theorem. At the moment, we do not know  examples of cardinality $2$ or $3$. \\
2. The property asserted in Theorem~\ref{1} (iv) implies that the natural morphism of pro-spaces
\[
X_{\et}(p) \longrightarrow K(G_S(p),1)
\]
from the pro-$p$-completion of the {\'e}tale homotopy type $X_{\et}$ of $X$ (see \cite{AM}) to the $K(\pi,1)$-pro-space attached to the pro-$p$-group $G_S(p)$ is a weak equivalence. Since $G_S(p)$ is the fundamental group of $X_{\et}(p)$, this justifies the notion  `$K(\pi,1)$ for $p$ and the {\'e}tale topology'. If $S$ contains the prime number $p$, this property always holds (cf.\ \cite{Sc2}).

\medskip
We can enlarge the set of prime numbers $S$ by the following

\begin{theorem} \label{up} Let $p$ be an odd prime number and let $S$ be a finite set of prime numbers congruent to $1$ modulo $p$. Assume that $\cd G_S(p)=2$. Let $\ell \notin S$ be another prime number congruent to $1$ modulo $p$ which does not split completely in the extension $\Q_S(p)/\Q$. Then $\cd G_{S \cup \{\ell\}}(p)=2$.
\end{theorem}

\section{Comparison with {\'e}tale cohomology}

In this section we show that cohomological dimension~$2$ implies the $K(\pi,1)$-property. 

\begin{lemma} \label{relrank}
Let $p$ be an odd prime number and let $S$ be a finite set of prime numbers congruent to $1$ modulo $p$.  Then
\[
\text{dim}_{{\mathbb F}_p} H^i(G_S(p),\Z/p\Z)=\left\{ 
\begin{array}{cl}
1& \text{if }\, i=0\\
\#S& \text{if }\, i=1\\
\#S& \text{if }\, i=2\, .
\end{array}
\right.
\]
\end{lemma}

\begin{proof}
The statement for $H^0$ is obvious.  \cite{NSW}, Theorem 8.7.11 implies the statement on $H^1$ and yields the inequality
\[
\text{dim}_{{\mathbb F}_p} H^2(G_S(p),\Z/p\Z) \leq \#S\, .
\] 
The abelian pro-$p$-group $G_S(p)^{\ab}$ has $\#S$ generators. There is only one $\Z_p$-extension of $\Q$, namely the cyclotomic $\Z_p$-extension, which is ramified at $p$. Since $p$ is not in $S$, $G_S(p)^{\ab}$ is finite, which implies that $G_S(p)$ must have at least as many relations as generators.  By \cite{NSW}, Corollary 3.9.5, the relation rank of $G_S(p)$ is $\text{dim}_{{\mathbb F}_p} H^2(G_S(p),\Z/p\Z)$, which yields the remaining inequality for $H^2$.
\end{proof}

\begin{proposition} \label{comp} Let $p$ be an odd prime number and let $S$ be a finite set of prime numbers congruent to $1$ modulo $p$.  If $\cd G_S(p)\leq 2$, then the scheme $X=\Spec(\Z)-S$ is a $K(\pi,1)$ for $p$ and the {\'e}tale topology, i.e.\ for any discrete $p$-primary $G_S(p)$-module $M$, considered as locally constant {\'e}tale sheaf on $X$, the natural homomorphism
\[
H^i (G_S(p), M) \rightarrow H^i_{\et} (X,M)
\]
is an isomorphism for all $i$.
\end{proposition}

\begin{proof} 
Let $L/k$ be a finite subextension of $k$ in $k_S(p)$. We denote the normalization of $X$ in $L$ by $X_L$. Then $H^i_{\et} (X_L,\Z/p\Z)=0$ for $i>3$ (\cite{Ma}, \S 3, Proposition C). Since flat and {\'e}tale cohomology coincide for finite {\'e}tale group schemes (\cite{Mi1}, III, Theorem~3.9), the flat duality theorem of Artin-Mazur (\cite{Mi2}, III Theorem~3.1) implies
\[
H^3_{\et}(X_L,\Z/p\Z)=H^3_{\fl}(X_L,\Z/p\Z)\cong H^0_{\fl,c}(X_L,\mu_p)^\vee=0,
\]
since a $p$-extension of $\Q$ cannot contain a primitive $p$-th root of unity.
Let $\tilde X$ be the universal (pro-)$p$-covering of $X$. We consider the Hochschild-Serre spectral sequence
\[
E_2^{pq}=H^p(G_S(p),H^q_{\et}(\tilde X,\Z/p\Z))\Rightarrow H^{p+q}_{\et}(X,\Z/p\Z).
\]
{\'E}tale cohomology commutes with inverse limits of schemes if the transition maps are affine (see \cite{AGV}, VII, 5.8). Therefore we have $H^i_{\et}(\tilde X,\Z/p\Z)=0$ for $i\geq 3$, and for $i=1$ by definition. Hence $E_2^{ij}=0$ unless $i=0,2$. Using the assumption $\cd G_S(p)\leq 2$,
the spectral sequence implies isomorphisms $H^i(G_S(p),\Z/p\Z) \stackrel{\sim}{\to} H^i_{\et}(X,\Z/p\Z)$ for $i=0,1$ and a short exact sequence
\[
0 \to H^2(G_S(p),\Z/p\Z) \stackrel{\phi}{\to} H^2_{\et}(X,\Z/p\Z)\to  H^2_{\et}(\tilde X,\Z/p\Z)^{G_S(p)} \to 0.
\]
Let $\bar X=\Spec(\Z)$. By the flat duality theorem of Artin-Mazur, we have an iso\-mor\-phism  $H^2_{\et}(\bar X,\Z/p\Z)\cong H^1_{\fl}(\bar X, \mu_p)^\vee$. The flat Kummer sequence $0\to \mu_p \to {\mathbb G}_m \to {\mathbb G}_m \to 0$, together with $H^0_\fl(\bar X,{\mathbb G}_m)/p=0=\text{}_p H^1_\fl(\bar X,{\mathbb G}_m)$ implies $H^2_{\et}(\bar X,\Z/p\Z)=0$. Furthermore, $H^3_{\et}(\bar X,\Z/p\Z)\cong H^0_{\fl}(\bar X, \mu_p)^\vee=0$. Considering the {\'e}tale excision sequence for the pair $(\bar X,X)$, we obtain an isomorphism
\[
H^2_{\et}(X,\Z/p\Z) \stackrel{\sim}{\longrightarrow} \bigoplus_{\ell \in S} H^3_{\ell}(\Spec(\Z_\ell), \Z/p\Z).
\]
The local duality theorem (\cite{Mi2}, II, Theorem 1.8) implies 
\[
H^3_{\ell}(\Spec(\Z_\ell), \Z/p\Z)\cong \text{Hom}_{\text{Spec}(\Z_\ell)}(\Z/p\Z,{\mathbb G}_m)^\vee.
\]
All primes $\ell \in S$ are congruent to $1$ modulo $p$ by assumption, hence $\Z_\ell$ contains a primitive $p$-th root of unity for $\ell \in S$, and we obtain $
\text{dim}_{{\mathbb F}_p} H^2_{\et}(X,\Z/p\Z)=\#S.
$ 
Now Lemma~\ref{relrank} implies that $\phi$ is an isomorphism. We therefore obtain 
\[
H^2_{\et}(\tilde X,\Z/p\Z)^{G_S(p)}=0.
\]
Since $G_S(p)$ is a pro-$p$-group, this implies (\cite{NSW}, Corollary 1.7.4) that 
\[
H^2_{\et}(\tilde X,\Z/p\Z)=0.
\]
We conclude that the Hochschild-Serre spectral sequence degenerates to a series of isomorphisms
\[
H^i(G_S(p),\Z/p\Z) \stackrel{\sim}{\longrightarrow} H^i_{\et}(X,\Z/p\Z), \qquad i\geq 0.
\]
If $M$ is a finite $p$-primary $G_S(p)$-module, it has a composition series with graded pieces isomorphic to $\Z/p\Z$ with trivial $G_S(p)$-action (\cite{NSW}, Corollary 1.7.4), and the statement of the proposition for $M$ follows from that for $\Z/p\Z$ and from the five-lemma. An arbitrary discrete $p$-primary $G_S(p)$-module is the filtered inductive limit of finite $p$-primary $G_S(p)$-modules, and the statement of the proposition follows since group cohomology (\cite{NSW}, Proposition 1.5.1) and {\'e}tale cohomology (\cite{AGV}, VII, 3.3) commute with filtered inductive limits.
\end{proof}

\section{Proof of Theorem~\ref{1}}
In this section we prove Theorem~\ref{1}. Let $p$ be an odd prime number and let $S$ be a finite set of prime numbers congruent to $1$ modulo $p$. Assume that $G_S(p)\neq 1$ and $\cd G_S(p)\leq 2$.

\medskip
Let $U \subset G_S(p)$ be an open subgroup. 
The abelianization $U^{\ab}$ of $U$ is a finitely generated abelian pro-$p$-group. If $U^{\ab}$ were infinite, it would have a quotient isomorphic to $\Z_p$, which corresponds to a $\Z_p$-extension $K_\infty$ of the number field $K=\Q_S(p)^U$ inside $\Q_S(p)$.  
By \cite{NSW}, Theorem 10.3.20 (ii), a $\Z_p$-extension of a number field is ramified at at least one prime dividing $p$. This contradicts $K_\infty\subset \Q_S(p)$ and we conclude that $U^{\ab}$ is finite. 

\medskip
In particular,  $G_S(p)^{\ab}$ is finite. Hence $G_S(p)$ is not free, and we obtain $\cd G_S(p)=2$. This shows the first part of assertion (i) of Theorem~\ref{1}  
and assertion (iv) follows from Proposition~\ref{comp}.

\medskip
By Lemma~\ref{relrank}, we know that for each prime number $\ell \in S$, the group $G_{S\backslash \{\ell\}}(p)$ is a proper quotient of $G_S(p)$, hence each $\ell\in S$ is ramified in the extension $\Q_S(p)/\Q$.  Let $G_\ell(\Q_S(p)/\Q)$ denote the decomposition group of $\ell$ in $G_S(p)$ with respect to some prolongation of $\ell$ to $\Q_S(p)$. As a subgroup of $G_S(p)$, $G_\ell(\Q_S(p)/\Q)$  has cohomological dimension less or equal to $2$. 
We have a natural surjection $G(\Q_\ell(p)/\Q_\ell)\twoheadrightarrow G_\ell(\Q_S(p)/\Q)$. By \cite{NSW}, Theorem 7.5.2, $G(\Q_\ell(p)/\Q_\ell)$ is the pro-$p$-group on two generators $\sigma, \tau$ subject to the relation $\sigma\tau\sigma^{-1}=\tau^\ell$. $\tau$ is a generator of the inertia group and $\sigma$ is a Frobenius lift. Therefore, $G(\Q_\ell(p)/\Q_\ell)$ has only three quotients of cohomological dimension less or equal to $2$: itself, the trivial group and the Galois group of the maximal unramified $p$-extension of $\Q_\ell$.  Since $\ell$ is ramified in the extension $\Q_S(p)/\Q$, the map $G(\Q_\ell(p)/\Q_\ell)\twoheadrightarrow G_\ell(\Q_S(p)/\Q)$ is an isomorphism, and hence $\Q_S(p)$ realizes the maximal $p$-extension of\/ $\Q_\ell$. This shows statement (iii) of Theorem~\ref{1}.

\medskip
Next we show the second part of statement (i). By \cite{NSW}, Proposition 3.3.3, we have $\scd G_S(p) \in \{2,3\}$. Assume that $\scd G=2$. We consider the $G_S(p)$-module
\[
D_2(\Z)=\varinjlim_{U} U^{\ab},
\]
where the limit runs over all open normal subgroups $U\lhd G_S(p)$ and for $V \subset U$ the transition map is the transfer $\text{Ver}\colon U^{\ab} \to V^{\ab}$, i.e.\ the dual of the corestriction map $\text{cor}\colon H^2(V,\Z) \to H^2(U,\Z)$ (see \cite{NSW}, I, \S5). 
By \cite{NSW}, Theorem 3.6.4 (iv), we obtain $G_S(p)^{\ab}=D_2(\Z)^{G_S(p)}$. On the other hand, $U^{\ab}$ is finite for all $U$ and the group theoretical version of the Principal Ideal Theorem (see \cite{Ne}, VI, Theorem 7.6) implies $D_2(\Z)=0$. Hence $G_S(p)^{\ab}=0$ which implies $G_S(p)=1$ producing a contradiction. Hence $\scd G_S(p)=3$ showing the remaining assertion of  Theorem~\ref{1}, (i). 

\medskip
It remains to show that $G_S(p)$ is a duality group. By \cite{NSW}, Theorem 3.4.6, it suffices to show that the terms
\[
D_i(G_S(p),\Z/p\Z)= \varinjlim_{U} H^i(U,\Z/p\Z)^\vee
\]
are trivial for $i=0,1$. Here $U$ runs through the open subgroups of $G_S(p)$, $^\vee$ denotes the Pontryagin dual and the transition maps are the duals of the corestriction maps. For $i=0$, and $V\subsetneqq U$, the transition map 
\[
\text{cor}^\vee\colon \Z/p\Z=H^0(V,\Z/p\Z)^\vee \to H^0(U,\Z/p\Z)^\vee=\Z/pZ
\]
is multiplication by $(U:V)$, hence zero. Since $G_S(p)$ is infinite, we obtain $D_0(G_S(p),\Z/p\Z)=0$. Furthermore,
\[
D_1(G_S(p),\Z/p\Z)=\varinjlim_U U^{\ab}/p=0
\]
by the Principal Ideal Theorem. This finishes the proof of Theorem~\ref{1}.

\section{The dualizing module}

Having seen that $G_S(p)$ is a duality group under certain conditions, it is interesting to calculate its dualizing module. The aim of this section is to prove

\begin{theorem} \label{dualthm} Let $p$ be an odd prime number and let $S$ be a finite set of prime numbers congruent to $1$ modulo $p$. Assume that $\cd G_S(p)=2$. Then we have a natural isomorphism  
\[
D\cong \text{\rm tor}_p\big(C_S(\Q_S(p))\big)
\]
between the dualizing module $D$ of\/ $G_S(p)$ and the $p$-torsion submodule of the $S$-id{\`e}le class group of\/ $\Q_S(p)$. There is  a natural short exact sequence
\[
0 \to \bigoplus_{\ell \in S} \text{\rm Ind\,}^{G_\ell}_{G_S(p)}\ \mu_{p^\infty} (\Q_\ell(p)) \to  D \to E_S(\Q_S(p)) \otimes \Q_p/\Z_p \to 0,
\]
in which $G_\ell$ is the decomposition group of\/ $\ell$ in $G_S(p)$ and $E_S(\Q_S(p))$ is the group of\/ $S$-units of the field $\Q_S(p)$.
\end{theorem}

\medskip
Working in a more general situation, let $S$ be a non-empty set of primes of a number field $k$.  We recall  some well-known facts from class field theory and we give some modifications for which we do not know a good reference.

\medskip\noindent
By $k_S$ we denote the maximal extension of $k$ which is unramified outside $S$ and we denote $G(k_S/k)$ by $G_S(k)$. For an intermediate field $k \subset K \subset k_S$, let $C_S(K)$ denote the $S$-id{\`e}le class group of $K$.
If $S$ contains the set $S_\infty$ of archimedean primes of $k$, then the pair $(G_S(k), C_S(k_S))$ is a class formation, see \cite{NSW}, Proposition~8.3.8. This remains true for arbitrary non-empty $S$, as can be seen as follows: We have the class formation 
\[
(G_S(k), C_{S \cup S_\infty}(k_{S})).
\]
Since $k_S$ is closed under unramified extensions, the Principal Ideal Theorem implies $\text{\it Cl}_S(k_S)=0$. Therefore we obtain the exact sequence
\[
0\to \bigoplus_{v\in S_\infty \backslash S (k)}\text{Ind}_{G_S(k)} k_v^\times \to C_{S \cup S_\infty}(k_{S}) \to C_S(k_S)\to 0.
\]
Since the left term is a cohomologically trivial $G_S(k)$-module, we obtain that $(G_S(k),C_S(k_S))$ is a class formation.
Finally, if $p$ is a prime number, then also $(G_S(k)(p),C_S(k_S(p))$ is a class formation.

\medskip\noindent
{\it Remark:} An advantage of considering the class formation $(G_S(k)(p),C_S(k_S(p))$ for sets $S$ of primes which do not contain $S_\infty$ is that we get rid of `redundancy at infinity'. A technical disadvantage is the absence of a reasonable Hausdorff topology on the groups $C_S(K)$ for finite subextensions $K$ of $k$ in $k_S(p)$.

\medskip

Next we calculate the module
\[
D_2(\Z_p)= \varinjlim_{U,n} H^2(U,\Z/p^n\Z)^\vee,
\]
where $n$ runs through all natural numbers, $U$ runs through all open subgroups of $G_S(k)(p)$ and $^\vee$ is the Pontryagin dual. If $\cd G_S(p)=2$, then $D_2(\Z_p)$ is the dualizing module $D$ of\/ $G_S(k)(p)$.

\begin{theorem} \label{kpicor} Let $k$ be a number field, $p$ an odd prime number and $S$ a finite non-empty set of non-archimedean primes of\/ $k$ such that the norm $N({\mathfrak p})$ of $\mathfrak p$ is congruent to $1$ modulo $p$ for all\/ ${\mathfrak p}\in S$. Assume that the scheme $X=\Spec({\cal O}_k)-S$ is a $K(\pi,1)$ for $p$ and the {\'e}tale topology and that $k_S(p)$ realizes the maximal $p$-extension $k_{\mathfrak p}(p)$ of\/ $k_{\mathfrak p}$ for all\/ ${\mathfrak p}\in S$. 
Then $G_S(p)$ is a pro-$p$-duality group of dimension $2$ with dualizing module
\[
D \cong \text{\rm tor}_{p}\big(C_S(k_S(p)\big).
\] 
\end{theorem}

\noindent
{\it Remarks.} 1. In view of Theorem~\ref{1}, Theorem~\ref{kpicor} shows Theorem~\ref{dualthm}.\\
2. In the case when $S$ contains all primes dividing $p$, a similar result has been proven in \cite{NSW}, X, \S5.

\begin{proof}[Proof of Theorem~\ref{kpicor}] 
We consider the schemes $\bar X=\Spec({\cal O}_k)$ and $X=\bar X -S$ and  we denote the natural embedding by $j: X \to \bar X$.  As in the proof of Proposition~\ref{comp}, the flat duality theorem of Artin-Mazur implies 
\[
H^3_{\et}(X,\Z/p\Z)\cong H^0_{\fl, c}(X,\mu_p)^\vee,
\] 
and the group on the right vanishes since $k_{\mathfrak p}$ contains a primitive $p$-th root of unity for all ${\mathfrak p}\in S$. The $K(\pi,1)$-property yields $\cd G_S(k)(p)\leq 2$. Since $k_S(p)$ realizes the maximal   $p$-extension $k_{\mathfrak p}(p)$ of $k_{\mathfrak p}$ for all ${\mathfrak p}\in S$, the inertia groups of these primes are of cohomological dimension~$2$ and we obtain $\cd G_S(p)=2$. 

Next we consider, for some $n\in {\mathbb N}$, the constant sheaf $\Z/p^n\Z$ on $X$.
The duality theorem of Artin-Verdier shows an isomorphism
\[
H^i_{\et}(\bar X,j_!(\Z/p^n\Z)) = H^i_c(X,\Z/p^n\Z)\cong {\text Ext}^{3-i}_X(\Z/p^n\Z,{\mathbb G}_m)^\vee.
\]
For ${\mathfrak p} \in S$, a standard calculation (see, e.g., \cite{Mi2}, II, Proposition 1.1) shows
\[
H^i_{\mathfrak p}(\bar X, j_!(\Z/p^n\Z)\cong H^{i-1}(k_{\mathfrak p},\Z/p^n\Z),
\]
where $k_{\mathfrak p}$ is (depending on the readers preference) the henselization or the completion of $k$ at $\mathfrak p$. 
The excision sequence for the pair $(\bar X,X)$ and the sheaf $j_!(\Z/p^n\Z)$ therefore implies a long exact sequence
\[
(\ast)\hspace{.2cm}  \cdots \to H^i_{\et}(X,\Z/p^n\Z) \to \bigoplus_{{\mathfrak p}\in S} H^i(k_{\mathfrak p},\Z/p^n\Z) \to
\text{Ext}_X^{2-i}(\Z/p^n\Z,{\mathbb G_m})^\vee \to \cdots
\]
The local duality theorem (\cite{NSW}, Theorem 7.2.6)  yields isomorphisms
\[
H^i(k_{\mathfrak p},\Z/p^n\Z)^\vee \cong H^{2-i}(k_{\mathfrak p},\mu_{p^n})
\]
for all $i\in \Z$.
Furthermore, 
\[
\text{Ext}_X^{0}(\Z/p^n\Z,{\mathbb G_m})=H^0(k,\mu_{p^n}).
\]
We denote by $E_S(k)$ and  $\text{\it Cl}_S(k)$ the group of $S$-units and the $S$-ideal class group of\/ $k$, respectively. By $\text{\it Br}(X)$, we denote the Brauer group of $X$.  The short exact sequence $0\to \Z \to \Z\to \Z/p^n\Z\to 0$ together with 
\[
\text{Ext}^i_X(\Z,{\mathbb G}_m)=H^i_{\et}(X,{\mathbb G}_m)=\left\{\begin{array}{cc}
E_S(k)& \text{for } i=0\\
\text{\it Cl}_S(k)&\text{for } i=1\\
\text{\it Br}(X) &\text{for } i=2
\end{array}
\right.
\]
and the Hasse principle for the Brauer group implies  exact sequences
\[
0 \to E_S(k)/p^n \to \text{Ext}_X^{1}(\Z/p^n\Z,{\mathbb G_m}) \to \text{ }_{p^n} \text{\it Cl}_S(k) \to 0
\]
and 
\[
0 \to \text{\it Cl}_S(k)/p^n \to \text{Ext}_X^{2}(\Z/p^n\Z,{\mathbb G_m}) \to \bigoplus_{{\mathfrak p}\in S} \text{ }_{p^n}\text{\it Br}(k_{\mathfrak p}).
\]
The same holds, if we replace $X$ by its normalization $X_K$ in a finite extension $K$ of $k$ in $k_S(p)$. Now we go to the limit over all such $K$.
Since $k_S(p)$ realizes the maximal $p$-extension of $k_{\mathfrak p}$ for all ${\mathfrak p}\in S$, we have 
\[
\varinjlim_{K} \bigoplus_{{\mathfrak p}\in S(K)} H^i(K_{\mathfrak p},\Z/p^n\Z)^\vee =
\varinjlim_{K} \bigoplus_{{\mathfrak p}\in S(K)} H^i(K_{\mathfrak p},\mu_{p^n})=0
\]
for $i\geq 1$ and
\[
\varinjlim_{K} \bigoplus_{{\mathfrak p}\in S(K)} \text{ }_{p^n}\text{\it Br}(K_{\mathfrak p})=0.
\]
The Principal Ideal Theorem implies $\text{\it Cl}_S(k_S(p))/p=0$ and since this group is a torsion group, its $p$-torsion part is trivial. Going to the limit over the exact sequences $(\ast)$ for all $X_K$, we obtain $D_i(\Z/p\Z)=0$ for $i=0,1$, hence $G_S(k)(p)$ is a duality group of dimension~$2$. Furthermore, we obtain the exact sequence

\bigskip
$\displaystyle 0 \to \text{tor}_p\big(E_S(k_S(p))\big) \to \bigoplus_{{\mathfrak p} \in S} \text{\rm Ind\,}^{G_{\mathfrak p}}_{G_S(k)(p)}\ \text{tor}_p\big(k_{\mathfrak p}(p)^\times\big) \to 
$
\begin{flushright}
$D \to E_S(k_S(p)) \otimes \Q_p/\Z_p \to 0.$
\end{flushright}

\medskip\noindent
Let $U\subset G_S(k)(p)$ be an open subgroup  and put $K=k_S(p)^U$. The invariant map
\[
\text{inv}_K\colon H^2(U,C_S(k_S(p))) \to \Q/\Z
\]
induces a pairing
\[
\text{Hom}_U(\Z/p^n\Z,C_S(k_S(p))) \times H^2(U,\Z/p^n\Z) \stackrel{\cup}{\to} H^2(U,C_S(K)) \stackrel{\text{inv}_K}{\to} \Q/\Z,
\]
and therefore a compatible system of maps 
\[
\text{}_{p^n}C_S(K) \to H^2(U,\Z/p^n\Z)^\vee
\]
for all $U$ and $n$. In the limit, we obtain a natural map
\[
\phi\colon \text{tor}_p\big(C_S(k_S(p)\big) \longrightarrow D.
\]
By our assumptions, the id{\`e}le group $J_S(k_S(p))$ is $p$-divisible. We therefore obtain the exact sequence

\bigskip
$\displaystyle 0 \to \text{tor}_p\big(E_S(k_S(p))\big) \to \bigoplus_{{\mathfrak p} \in S} \text{\rm Ind\,}^{G_{\mathfrak p}}_{G_S(k)(p)}\ \text{tor}_p\big(k_{\mathfrak p}(p)^\times\big) \to 
$
\begin{flushright}
$ \text{tor}_p\big(C_S(k_S(p))\big) \to E_S(k_S(p)) \otimes \Q_p/\Z_p \to 0$
\end{flushright}

\medskip\noindent
which, via the just constructed map $\phi$, compares to the similar sequence with $D$ above. Hence $\phi$ is an isomorphism by the five lemma.
\end{proof}

Finally, without any assumptions on $G_S(k)(p)$, we calculate the  $G_S(k)(p)$-module $D_2(\Z_p)$  as a quotient of $\text{tor}_p \big(C_S(k_S(p))\big)$ by a subgroup of universal norms. We therefore can interpret Theorem~\ref{kpicor} as a vanishing statement on universal norms.

\medskip
Let us fix some notation.  
If $G$ is a profinite group and if $M$ is a $G$-module, we denote by $_{p^n}M$ the submodule of elements annihilated by $p^n$. By $N_G(M)\subset M^G$ we denote the subgroup of universal norms, i.e.
\[
N_G(M)= \bigcap_{U}\ N_{G/U} (M^U),
\]
where $U$ runs through the open normal subgroups of $G$ and $N_{G/U}(M^U) \subset M^G$
is the image of the norm map
\[
N\colon M^U \to M^G, \ m\mapsto \sum_{\sigma \in G/U} \sigma m.
\]

\begin{proposition}\label{dual}
Let $S$ be a non-empty finite set of non-archimedean primes of\/ $k$ and let $p$ be an odd prime number such that $S$ contains no prime dividing~$p$. Then
\[
D_2(G_S(k)(p),\Z_p)\cong \varinjlim_{K,n} \text{ }_{p^n}C_S(K)/N_{G(k_S(p)/K)}(_{p^n}C_S(K)),
\]
where $n$ runs through all natural numbers and $K$ runs through all finite subextension of $k$ in $k_S(p)$.
\end{proposition}

\begin{proof} 
We want to use Poitou's duality theorem (\cite{Sc2}, Theorem~1). But the class module $C_S(k_S(p))$ is not level-compact and we cannot apply the theorem directly. Instead, we consider the level-compact class formation 
\[
\big(G_S(k)(p),C_{S\cup S_\infty}^0(k_S(p))\big),
\]
 where $C_{S\cup S_\infty}^0(k_S(p)) \subset C_{S\cup S_\infty}(k_S(p))$ is the subgroup of id{\`e}le classes of norm~$1$. 
By \cite{Sc2}, Theorem 1, we have for all natural numbers $n$ and all finite subextensions $K$ of $k$ in $k_S(p)$ a natural isomorphism
\[
H^2(G_S(K)(p),\Z/p^n\Z)^\vee \cong \hat H^0(G_S(K)(p),\  _{p^n}C_{S\cup S_\infty}^0(k_S(p))),
\]
where $\hat H^0$ is Tate-cohomology in dimension~$0$ (cf.\ \cite{Sc2}). The exact sequence
\[
0\to \bigoplus_{v\in S_\infty  (K)} K_v^\times \to C_{S \cup S_\infty}(K) \to C_S(K)\to 0
\]
and the fact that $K_v^\times$ is $p$-divisible for archimedean $v$, implies for all $n$ and all finite subextensions $K$ of $k$ in $k_S(p)$ an exact sequence of finite abelian groups 
\[
0\to \bigoplus_{v\in S_\infty (K)} \mu_{p^n}(K_v) \to \text{}_{p^n}C_{S \cup S_\infty}(K) \to \text{ }_{p^n}C_S(K)\to 0.
\]
\cite{Sc2}, Proposition~$7$ therefore implies  isomorphisms
\[
\hat H^0(G_S(K)(p),\text{}_{p^n}C_{S \cup S_\infty}(k_S(p))) \cong
\hat H^0(G_S(K)(p),\text{}_{p^n}C_{S}(k_S(p))) 
\]
for all $n$ and $K$. Furthermore, the exact sequence 
\[
0 \to C_{S\cup S_\infty}^0(K) \to C_{S\cup S_\infty}(K) \stackrel{|\;|}{\to} {\mathbb R}_+^\times \to
0
\]
shows $\text{}_{p^n}C^0_{S \cup S_\infty}(K))= \text{}_{p^n}C_{S \cup S_\infty}(K))$ for all $n$ and all finite subextensions $K$ of $k$ in $k_S(p)$. 
Finally,  \cite{Sc2}, Lemma~5 yields isomorphisms
\[
\hat H^0(G_S(K)(p),\text{}_{p^n}C_{S}(k_S(p))) \cong \text{ }_{p^n}C_S(K)/N_{G(k_S(p)/K)}(_{p^n}C_S(K)).
\]
Going to the limit over all $n$ and $K$, we obtain the statement of the Proposition.
\end{proof}

\pagebreak

\section{Going up}
The aim of this section is to prove Theorem~\ref{up}. 
We start with the following lemma.

\begin{lemma} \label{loccalc} Let $\ell\neq p$ be prime numbers. Let $\Q_\ell^h$ be the henselization of $\Q$ at $\ell$ and let $K$ be an algebraic extension of $\Q_\ell^h$ containing the maximal unramified $p$-extension $(\Q_\ell^h)^{nr,p}$ of\/ $\Q_\ell^h$. Let $Y=\Spec({\cal O}_K)$, and denote the closed point of\/ $Y$ by $y$.  Then the local {\'e}tale cohomology group $H^i_y(Y,\Z/p\Z)$ vanishes for $i\neq 2$ and we have a natural isomorphism
\[
H^2_y(Y,\Z/p\Z)\cong H^1(G(K(p)/K),\Z/p\Z).
\]
\end{lemma}

\begin{proof} Since $K$ contains $(\Q_\ell^h)^{nr,p}$, we have $H^i_{\et}(Y,\Z/p\Z)=0$ for $i >0$. The excision sequence shows $H^i_y(Y,\Z/p\Z)=0$ for $i=0,1$ and $H^i_y(Y,\Z/p\Z)\cong H^{i-1}(G(\bar K/K),\Z/p\Z)$  for $i\geq 2$. By \cite{NSW}, Proposition 7.5.7, we have
\[
H^{i-1}(G(\bar K/K),\Z/p\Z)=H^{i-1}(G(K(p)/K),\Z/p\Z)
\] 
But $G(K(p)/K)$ is a free pro-$p$-group (either trivial or isomorphic to $\Z_p$). This concludes the proof.
\end{proof}

Let $k$ be a number field and let $S$ be finite set of primes of $k$.  For a (possibly infinite) algebraic extension $K$ of $k$ we denote by $S(K)$ the set of prolongations of primes in $S$ to $K$. Now assume that $M/K/k$ is a tower of pro-$p$ Galois extensions.  We denote the inertia group of a prime ${\mathfrak p}\in S(K)$ in the extension $M/K$ by $T_{\mathfrak p}(M/K)$. For $i\geq 0$ we write
\[
\ressum_{{\mathfrak p} \in S(K)} H^i(T_{\mathfrak p}(M/K), \Z/p\Z)
\stackrel{df}{=} 
\varinjlim_{k'\subset K} \bigoplus_{{\mathfrak p} \in S(k')} H^i(T_{\mathfrak p}(M/k'), \Z/p\Z),
\]
where the limit on the right hand side runs through all finite subextensions $k'$ of $k$ in $K$. The $G(K/k)$-module $\ressumsmall_{{\mathfrak p} \in S(K)} H^i(T_{\mathfrak p}(M/K), \Z/p\Z)$ is the maximal discrete submodule of the  product $\prod_{{\mathfrak p} \in S(K)} H^i(T_{\mathfrak p}(M/K), \Z/p\Z) $.

\begin{proposition}\label{cohoben}
Let $p$ be an odd prime number and let $S$ be a finite set of prime numbers congruent to $1$ modulo $p$ such that $\cd G_S(p)=2$. Let $\ell \notin S$ be another prime number congruent to $1$ modulo $p$ which does not split completely in the extension $\Q_S(p)/\Q$. Then, for any prime $\mathfrak p$ dividing $\ell$ in $\Q_S(p)$, the inertia group of\/ $\mathfrak p$ in the extension $\Q_{S \cup \{\ell\}}(p)/\Q_S(p)$ is infinite cyclic. Furthermore, 
\[H^i(G(\Q_{S \cup \{\ell\}}(p)/\Q_S(p)),\Z/p\Z)=0
\] 
for $i\geq 2$. For $i=1$ we have a natural isomorphism
\[
H^1(G(\Q_{S \cup \{\ell\}}(p)/\Q_S(p)),\Z/p\Z) \cong \ressum_{{\mathfrak p} \in S_\ell(\Q_S(p))} \!\!\!\! H^1(T_{\mathfrak p}(\Q_{S \cup \{\ell\}}(p))/\Q_S(p),\Z/p\Z),
\]
where $S_\ell(\Q_S(p))$ denotes the set of primes of\/ $\Q_S(p)$ dividing $\ell$.
In particular, $G(\Q_{S \cup \{\ell\}}(p)/\Q_S(p))$ is a free pro-$p$-group.
\end{proposition}

\begin{proof} Since $\ell$ does not split completely in $\Q_S(p)/\Q$ and since $\cd G_S(p)=2$, the decomposition group of $\ell$ in $\Q_S(p)/\Q$ is a non-trivial and torsion-free quotient of $\Z_p\cong G(\Q_\ell^{nr,p}/\Q_\ell)$. Therefore $\Q_S(p)$ realizes the maximal unramified $p$-extension of $\Q_\ell$.
We consider the scheme $X=\Spec(\Z)-S$ and its universal pro-$p$ covering $\tilde X$ whose field of functions is $\Q_S(p)$. Let $Y$ be the subscheme of $\tilde X$ obtained by removing all primes of residue characteristic $\ell$. We consider the {\'e}tale excision sequence for 
the pair $(\tilde X, Y)$. By Theorem~\ref{comp}, we have $H^i_{\et}(\tilde X,\Z/p\Z)=0$ for $i >0$, which implies isomorphisms
\[
H^i_{\et}(Y,\Z/p\Z) \stackrel{\sim}{\to} \ressum_{{\mathfrak p} \mid \ell} H^{i+1}_{\mathfrak p}(Y_{\mathfrak p}^h,\Z/p\Z)
\]
for $i\geq 1$. By Lemma~\ref{loccalc}, we obtain $H^i_{\et}(Y,\Z/p\Z)=0$ for $i\geq 2$. The universal $p$-covering $\tilde Y$ of\/ $Y$ has $\Q_{S\cup \{\ell\}}(p)$ as its function field, and the Hochschild-Serre spectral sequence for $\tilde Y/Y$ yields an inclusion
\[
H^2(G (\Q_{S\cup \{\ell\}}(p)/\Q_{S}(p)),\Z/p\Z) \hookrightarrow H^2_{\et}(Y,\Z/p\Z)=0.
\]
Hence $G (\Q_{S\cup \{\ell\}}(p)/\Q_{S}(p))$ is a free pro-$p$-group and for $H^1$ we obtain

\medskip\noindent
$
\ H^1(G (\Q_{S\cup \{\ell\}}(p)/\Q_{S}(p)),\Z/p\Z)\stackrel{\sim}{\to} H^1_{\et}(Y,\Z/p\Z)
$
\begin{flushright}
$
\cong \ressum_{{\mathfrak p}\in S_\ell(\Q_S(p))} H^1(G(\Q_S(p)_{\mathfrak p}(p)/\Q_S(p)_{\mathfrak p}),\Z/p\Z).
$
\end{flushright}
This shows that each ${\mathfrak p} \mid \ell$ ramifies in $\Q_{S\cup \{\ell\}}(p)/\Q_{S}(p)$, and since the Galois group is free, $\Q_{S\cup \{\ell\}}(p)$ realizes the maximal $p$-extension of $\Q_S(p)_{\mathfrak p}$. In particular, 
\[
H^1(G(\Q_S(p)_{\mathfrak p}(p)/\Q_S(p)_{\mathfrak p}),\Z/p\Z) \cong H^1(T_{\mathfrak p}(\Q_{S \cup \{\ell\}}(p)/\Q_S(p)),\Z/p\Z)
\]
for all $\mathfrak p \mid \ell$, which finishes the proof.
\end{proof}

Let us mention in passing that the above calculations imply the validity of the following arithmetic form of Riemann's existence theorem.

\begin{theorem}
Let $p$ be an odd prime number and let $S$ be a finite set of prime numbers congruent to $1$ modulo $p$ such that $\cd G_S(p)=2$. Let $T\supset S$ be another set of prime numbers congruent to $1$ modulo $p$. Assume that all $\ell \in T\backslash S$ do not split completely in the extension $\Q_S(p)/\Q$. Then the inertia groups in $\Q_{T}(p)/\Q_S(p)$ of all primes $\mathfrak p \in T\backslash S(\Q_S(p))$ are infinite cyclic  and the natural homomorphism
\[
\phi:\freeproductmed_{{\mathfrak p} \in T\backslash S(\Q_S(p))} T_{\mathfrak p}(\Q_T(p)/\Q_S(p)) \longrightarrow G(\Q_{T}(p)/\Q_S(p))
\]
is an isomorphism.
\end{theorem}

\medskip\noindent
{\it Remark:} A similar theorem holds in the case that $S$ contains~$p$, see \cite{NSW}, Theorem~10.5.1.

\begin{proof} By Proposition~\ref{cohoben} and by the calculation of the cohomology of a free product (\cite{NSW}, 4.3.10 and 4.1.4), $\phi$ is a homomorphism between free pro-$p$-groups which induces an isomorphism on mod~$p$ cohomology. Therefore $\phi$ is an isomorphism.
\end{proof}

\begin{proof}[Proof of theorem~\ref{up}] We consider the Hochschild-Serre spectral sequence
\[
E_2^{ij}=H^i(G_{S}(p),H^j(G(\Q_{S\cup\{\ell\}}(p)/\Q_S(p)),\Z/p\Z) \Rightarrow H^{i+j}(G_{S\cup\{\ell\}}(p),\Z/p\Z).
\]
By Proposition~\ref{cohoben}, we have $E_2^{ij}=0$ for $j\geq 2$ and 
\[
H^1(G(\Q_{S\cup\{\ell\}}(p)/\Q_S(p)),\Z/p\Z)\cong \ressum_{{\mathfrak p}\mid \ell} 
H^1(T_{\mathfrak p}(\Q_{S \cup \{\ell\}}(p)/\Q_S(p)),\Z/p\Z)
\]
\[
\cong \text{Ind}^{G_\ell}_{G_S(p)} H^1(T_{\ell}(\Q_{S \cup \{\ell\}}(p)/\Q_S(p)),\Z/p\Z),
\]
where $G_\ell\cong \Z_p$ is the decomposition group of $\ell$ in $G_S(p)$. We obtain
$E^{i,1}_2=0$ for $i\geq 2$. By assumption, $\cd G_S(p)=2$, hence $E_2^{0,j}=0$ for $j \geq 3$. This implies $H^3(G_{S\cup\{\ell\}}(p),\Z/p\Z)=0$, and hence $\cd G_{S\cup\{\ell\}}(p)\leq 2$. Finally, the decomposition group of $\ell$ in $G_{S\cup\{\ell\}}(p)$ is full, i.e.\ of cohomological dimension~$2$. Therefore, $\cd G_{S\cup\{\ell\}}(p)= 2$.
\end{proof}

We obtain the following

\begin{corollary}\label{upcor}
Let $p$ be an odd prime number and let $S$ be a finite set of prime numbers congruent to $1$ modulo $p$. Let $\ell\notin S$ be a another prime number congruent to $1$ modulo $p$. Assume that there exists a prime number $q \in S$ such that the order of\/ $\ell$ in $(\Z/q\Z)^\times$ is divisible by $p$ (e.g.\ $\ell$ is not a $p$-th power modulo~$q$). Then $\cd G_S(p)=2$ implies $\cd G_{S\cup \{\ell\}}(p)=2$.
\end{corollary}

\begin{proof}
Let $K_q$ be the maximal subextension of $p$-power degree in $\Q(\mu_q)/\Q$. Then $K_q$ is a non-trivial finite subextension of\/ $\Q$ in $\Q_S(p)$ and $\ell$ does not split completely in $K_q/\Q$. Hence the result follows from Theorem~\ref{up}.
\end{proof}

\noindent
{\it Remark.} One can sharpen Corollary~\ref{upcor} by finding weaker conditions on a prime $\ell$ not to split completely in $\Q_S(p)$.

\section{Proof of Theorem~\ref{2}}\label{lab}

In this section we prove Theorem~\ref{2}.  We start by recalling the notion of the linking diagram attached to $S$ and $p$ from \cite{La}. Let $p$ be an odd prime number and let $S$ be a finite set of prime numbers congruent to $1$ modulo $p$. Let $\Gamma (S)(p)$ be the directed graph with vertices the primes of $S$ and edges the pairs $(r, s) \in  S \times S$ with $r$ not a $p$-th power modulo $s$. We now define a function $\ell$ on the set of pairs of distinct primes of $S$ with values in $\Z/p\Z$ by first choosing a primitive root $g_s$ modulo $s$ for each $s\in S$. Let $\ell_{rs} = \ell(r, s)$ be the image in $\Z/p\Z$ of any integer $c$ satisfying 
\[
r\equiv g_s^{-c} \bmod s\;.
\]
The residue class $\ell_{rs}$ is well-defined since $c$ is unique modulo $s-1$ and $p \mid s-1$. Note that $(r, s)$ is an edge of $\Gamma(S)(p)$ if and only if $\ell_{rs} \neq  0$. We call $\ell_{rs}$ the {\em linking number} of the pair $(r, s)$. This number depends on the choice of primitive roots, if $g$ is another primitive root modulo $s$ 
and $g_s \equiv g^a \bmod s$, then the linking number attached to $(r, s)$ would be multiplied by $a$ if $g$ were used instead of $g_s$. The directed graph $\Gamma(S)(p) $ together with $\ell$ is called the {\it linking diagram} attached to $S$ and $p$. 

\begin{definition} \rm
We call a finite set $S$ of prime numbers congruent to $1$ modulo $p$  {\em strictly circular with respect to $p$} (and $\Gamma(S)(p)$ a {\em non-singular circuit}), if there exists an ordering $S=\{q_1,\ldots,q_n\}$ of the primes in $S$ such that the following conditions hold.
\begin{enumerate}
\item[\rm (a)] The vertices $q_1,\ldots,q_n$ of $\Gamma(S)(p)$ form a circuit $q_1q_2\cdots q_nq_1$.
\item[\rm (b)] If $i,j$ are both odd, then $q_iq_j$ is not an edge of $\Gamma(S)(p)$.
\item[\rm (c)] If we put $\ell_{ij}=\ell(q_i,q_j)$, then
\[
\ell_{12}\ell_{23}\cdots \ell_{n-1,n}\ell_{n1} \neq 
\ell_{1n}\ell_{21}\cdots \ell_{n,n-1}.
\] 
\end{enumerate}
\end{definition}

\noindent
Note that condition (b) implies that $n$ is even $\geq 4$ and that (c) is satisfied if there is an edge $q_iq_j$ of the circuit $q_1q_2\cdots q_nq_1$ such that $q_jq_i$ is not an edge of $\Gamma (S)(p)$. Condition (c) is independent of the choice of primitive roots since the condition can be written in the form
\[ 
\frac{\ell_{1n}}{\ell_{n-1,n}} \frac{\ell_{21}}{\ell_{n1}} \frac{\ell_{32}}{\ell_{12}} \cdots \frac{\ell_{n,n-1}}{\ell_{n-2,n-1}} \neq 1,
\]
where each ratio in the product is independent of the choice of primitive roots.

\medskip
If $p$ is an odd prime number and if $S=\{q_1,\ldots,q_n\}$ is a finite set of prime numbers congruent to $1$ modulo $p$, then, by a result of Koch \cite{Ko}, the group $G_S(p)$ has a minimal presentation $G_S(p)=F/R$, where $F$ is a free pro-$p$-group on generators $x_1,\ldots,x_n$ and $R$ is the minimal normal subgroup in $F$ on generators $r_1,\ldots,r_n$, where
\[
r_i\equiv x_i^{q_i-1} \prod_{j\neq i} [x_i,x_j]^{\ell_{ij}} \bmod F_3.
\] 
Here $F_3$ is the third step of the lower $p$-central series of $F$ and the $\ell_{ij}=\ell(q_i,q_j)$ are the linking numbers for some choice of primitive roots. If $S$ is strictly circular, Labute (\cite{La}, Theorem~1.6) shows that $G_S(p)$ is a so-called `mild' pro-$p$-group, and, in particular, is of cohomological dimension~$2$ (\cite{La}, Theorem~1.2). 

\begin{proof}[Proof of Theorem~\ref{2}]
By \cite{La}, Theorem~1.6, we have $\cd G_T(p)=2$. By assumption, we find a series of subsets
\[
T=T_0\subset T_1 \subset \cdots \subset T_r=S,
\]
such that for all $i\geq 1$, the set $T_i\backslash T_{i-1}$ consists of a single prime number $q$ congruent to $1$ modulo $p$ and there exists a prime number $q' \in T_{i-1}$ with $q$ not a $p$-th power modulo $q'$. An inductive application of Corollary~\ref{upcor} yields the result.
\end{proof}

\noindent
{\it Remark.} Labute also proved some variants of his group theoretic result \cite{La}, Theorem~1.6.  The same proof as above shows corresponding variants of Theorem~\ref{2}, by replacing condition (i) by other conditions on the subset $T$ as they are described in \cite{La}, \S 3.

\bigskip
A straightforward applications of \v{C}ebotarev's density theorem shows that, given $\Gamma(S)(p)$, a prime number $q$ congruent to $1$ modulo~$p$ can be found with the additional edges of $\Gamma({S\cup \{q\}})(p)$ arbitrarily prescribed (cf.\ \cite{La}, Proposition~6.1). We therefore obtain the following corollaries.

\begin{corollary} Let $p$ be an odd prime number and let $S$ be a finite set of prime numbers congruent to $1$ modulo $p$, containing a strictly circular subset $T\subset S$. Then there exists a prime number $q$ congruent to $1$ modulo $p$ with
\[
\cd G_{S \cup \{q\}}(p)=2.
\]
\end{corollary}

\begin{corollary}
Let $p$ be an odd prime number and let $S$ be a finite set of prime numbers congruent to $1$ modulo $p$. Then we find a finite set $T$ of prime numbers congruent to $1$ modulo~$p$ such that
\[
\cd G_{S \cup T}(p)=2.
\]
\end{corollary}

\vskip1cm
\noindent \footnotesize{Alexander Schmidt, NWF I - Mathematik, Universit\"{a}t Regensburg, D-93040
Regensburg, Deutschland. email: alexander.schmidt@mathematik.uni-regensburg.de}
\end{document}